\documentclass [11pt] {amsart}

\newtheorem{theorem}{Theorem}[section]
\newtheorem{lemma}[theorem]{Lemma}

\newtheorem{remark}[theorem]{Remark}
\newtheorem{example}[theorem]{Example}

\numberwithin{equation}{section}

\newcommand{\C}{\mathbb{C}}

\newcommand{\iD}{\mathit{\Delta}}
\newcommand{\cE}{\mathcal{E}}
\newcommand{\cF}{\mathcal{F}}
\newcommand{\N}{\mathbb{N}}
\newcommand{\iO}{\mathit{\Omega}}
\newcommand{\iPhi}{\mathit{\Phi}}
\newcommand{\cP}{\mathcal{P}}
\newcommand{\R}{\mathbb{R}}
\newcommand{\cR}{\mathcal{R}}
\newcommand{\cS}{\mathcal{S}}
\newcommand{\iS}{\mathit{\Sigma}}
\newcommand{\supp}{\mathrm{supp}\,}
\newcommand{\cX}{\mathcal{X}}

\begin{document}

\title[Rational Approximation] 
      {Rational Approximation \\ for a Quasilinear Parabolic Equation}
\date{\today}


\author{P. M. Gauthier}
\address{D\'{e}partement de math\'{e}matiques et de statistique
\\
         Universit\'{e} de Montr\'{e}al
\\
         CP 6128 Centre Ville
\\
         Montr\'{e}al, Qu\'{e}bec H3C 3J7
\\
         Canada}

\email{gauthier@dms.umontreal.ca}


\author{N. Tarkhanov}
\address{Universit\"{a}t Potsdam
\\
         Institut f\"{u}r Mathematik
\\
         Am Neuen Palais 10
\\
         14469 Potsdam
\\
         Germany}

\email{tarkhanov@math.uni-potsdam.de}



\begin{abstract}
Approximation theorems, analogous to known results for linear elliptic
equations, are obtained for solutions of the heat equation.
Via the Cole-Hopf transformation, this gives rise to approximation theorems
for a nonlinear parabolic equation, Burgers' equation.
\end{abstract}

\maketitle

\tableofcontents

\section{Introduction}
\label{s.Int}

It is well known that each holomorphic function on a compact set $K$ in the
complex plane can be approximated by rational fractions with poles away from
$K$.
In fact the set of poles can be fixed arbitrarily to satisfy the only
condition that it meets each connected component of the complement of $K$.
If it has at least one limit point at each component of $\C \setminus K$,
then it is actually possible to approximate by finite linear combinations of
the Cauchy kernel.

The Cauchy kernel is a fundamental solution for the Cauchy-Riemann operator in
$\C$.
This gives insight into what can be though of as rational solutions to linear
partial differential equations.
These are just finite linear combinations of a right fundamental solution to
the given equation.
In this way the rational approximation theory naturally extends to solutions of
linear elliptic equations, see \cite[5.3.2]{T}.

Since the analogue of rational approximation has been well developed for
elliptic equations, it is natural, from the theoretical viewpoint, to attempt
a similar theory for parabolic equations.
Thus, we seek to approximate a solution of a parabolic equation by solutions
having isolated singularities, and we also attempt to specify the location of
these singularities.

As a possible physical interpretation for rational parabolic approximation,
consider the following.
Let $B$ be an open set in $\R^n$ and
let $u$ be a heat distribution on $B$, that is a solution of the heat equation
on $\R \times B \subset \R^{n+1}$.
We call the pair $(B,u)$ a thermal box.
Thus, a thermal box is a domain endowed with a heat distribution.
We may think of ovens and refrigerators as examples, but we may also think of
a house or a region of physical space as examples.
In the case that $B$ is a refrigerator, an oven or a house, we are usually
interested in prescribing a temperature distribution $u (t,x)$, which is
constant as a function of $x$.
But we could ask for much more.
We could ask for different temperatures in different parts of $B$ and, indeed,
this is often done in houses.
For example, let $B_1$ and
                 $B_2$ be two subregions of $B$ and let $u_1$ and
                                                        $u_2$
be heat distributions on $B_1$ and
                         $B_2$,
respectively.
Is it possible to design the thermal box $(B,u)$ to these specifications?
That is, can we find a temperature distribution $u$ on $B$ such that
   $u = u_1$ on $B_1$ and
   $u = u_2$ on $B_2$?
Since the heat distribution $u$ is analytic in the space variable, the answer
is in general ``no.''
However, from the point of view of every conceivable physical application, it
would be sufficient to find $u_{\varepsilon}$ on $B$, which approximate $u_1$
                                                                    and $u_2$
on $B_1$ and
   $B_2$,
respectively, within some tolerance $\varepsilon$, provided we can do this for
arbitrarily small tolerances $\varepsilon$.
This is precisely the kind of approximation we wish to investigate:
uniform approximation on subsets of $B$ by solutions of the heat equation on
all of $B$.
From the engineering point of view, one is interested
   not merely in the mathematical existence of the global approximating heat
   distribution $u_{\varepsilon}$.
One seeks a way of designing $B$ to obtain such a heat distribution
$u_{\varepsilon}$.
A mathematical model for this is to take as $u$ a finite linear combination of
fundamental solutions whose singularities lie on the boundary of $B$.
We shall call these simple rational functions.
We may think of such a boundary point as a heat source, or infinitely hot spot,
if the coefficient is positive at that point, and we may think of it as a heat
sink or infinitely cold spot, if the coefficient is negative.
As an engineering approximation, one can impose hot temperatures at the hot
spots and cold temperatures at the cold spots on the walls of the thermal box
$B$.

In addition to proving results on approximation by rational solutions to
parabolic equations, we shall also investigate the problem of approximating
solutions in a one domain by solutions in a given larger domain.
Pairs of domains where this is possible are called Runge pairs.

In 1929, Birkhoff \cite{B} showed the existence of an entire function $u (z)$
with the following remarkable property.
For any entire function $f (z)$, there is a sequence $\{ a_j \}$ such that the
sequence $\{ u (z + a_j) \}$ of translates of $u$ converges to the function
$f (z)$ uniformly on compact subsets of $\C$.
Thus, the translates of the function $u$ approximate all entire functions
$f$.
Such a function $u$ is said to be universal.

Another kind of universality is that of universal series.
Let $V$ be a space of functions and
    $\{ u_j \}$ a sequence of functions in $V$.
The series $\sum u_j$ is said to be universal in $V$ if the set formed by the
partial sums is dense in $V$.
Of course, such series are to be thought of as formal series.
They are never convergent, if $V$ is not trivial.
In 1951, Seleznev \cite{S} showed the existence of a universal power series
and this result was refined in \cite{L},
                               \cite{CP} and
                               \cite{N}.
Subsequently, other universal series were found for solutions to
   the Cauchy-Riemann as well as
   the Laplace equation
and other elliptic equations.
Recently, an abstract theory of universal series,
   covering most of these results,
were introduced in \cite{NP}.
For a broader view of universality, we refer the reader to the excellent
survey \cite{G}.

This work is intended as an attempt to specify these results in the context of
nonlinear partial differential equations.
We shall present results on universal functions as well as universal (formal)
series for solutions to parabolic equations.
The crucial role in our investigation is due to the so-called Cole-Hopf
transformation.

\section{The Cole-Hopf transformation}
\label{s.C-Ht}

The Cole-Hopf transformation was discovered independently by Cole \cite{C}
                                                         and Hopf \cite{H}
around 1950.
It changes Burgers' equation $u_t + u u_x = u_{xx}$ into the heat equation
$v_t = v_{xx}$.
To derive the transform, we let $u = - p_x$.
Then Burgers' equation can be integrated yielding
   $p_t - p_x^2 / 2 = p_{xx}$
up to a function depending on $t$ only.
Let $p = 2 \log v$.
Applying some algebra to this we get $v_t = v_{xx}$.

The $n\,$-dimensional forced Burgers equation
$
   u_t + (u \cdot \nabla) u = \nu \nabla^2 u - \nabla f (t,x)
$
for $u = - \nabla p$, which describes the dynamics of a stirred, pressure less
and vorticity-free fluid, has found interesting applications in a wide range
of non-equilibrium statistical physics problems, see \cite{BKh}.
Here, $\nu$ stands for the viscosity.
The associated Hamilton-Jacobi equation, satisfied by the velocity potential
$p$,
\begin{equation}
\label{eq.KPZ}
   p_t - \frac{1}{2} |\nabla p|^2 = \nu \nabla^2 p + f (t,x),
\end{equation}
has been frequently studied as a nonlinear model for the motion of an interface
under deposition, when the forcing potential $f$ is random, delta-correlated
in both space and time.
Equation (\ref{eq.KPZ}) is the well-known Kardar-Parisi-Zhang equation, see
   \cite{KPZh}.

Starting with this example, we consider a quasilinear partial differential
equation
\begin{equation}
\label{eq.qlpde}
   p_t = \iD p + a (p) |\nabla p|^2
\end{equation}
in $\R^{n+1}$, where $a$ is a continuous real-valued function on the real axis.
Choose a strictly monotone increasing $C^2$ function $u = U (p)$ on $\R$ with
the property that
$$
   a (p) = \frac{U'' (p)}{U' (p)}
$$
for all $p \in \R$.

The general solution of this ordinary differential equation satisfying the
initial condition $U' (0) = U_1 > 0$ is
\begin{equation}
\label{eq.symmetry}
   U' (p) = U_1\, \exp \Big( - \int_{0}^{p} a (t) dt \Big),
\end{equation}
which is a smooth function on $\R$ with positive values.
The function $u = U (p)$ may be found by integration.
In this way we recover what is referred to as the Cole-Hopf transformation.

A simple computation shows that the change of variables $u = U (p)$ reduces
(\ref{eq.qlpde}) to the heat equation
\begin{equation}
\label{eq.heat}
   u_t = \iD u
\end{equation}
for the new unknown function $u$.
Hence, the general solution to (\ref{eq.qlpde}) is $p (x) = U^{-1} (u (x))$,
with $u$ satisfying (\ref{eq.heat}).

\begin{example}
\label{e.CHT}
{\em
Let $a$ be constant.
Then
\begin{eqnarray*}
   U (p)
 & = &
   U_0 + U_1\, \frac{1 - \exp (- a p)}{a},
\\
   U^{-1} (u)
 & = &
 - \frac{1}{a} \log \Big( 1 - a\, \frac{u-U_0}{U_1} \Big).
\end{eqnarray*}
}
\end{example}

Using the function $U$ allows one to endow the set of solutions to equation
(\ref{eq.qlpde}) with the symmetry
$
   p_1 \circ p_2 := U^{-1} (U (p_1) + U (p_2)).
$

\section{Rational approximation}
\label{s.ra}

If $U$ is an open subset of $\R^{n+1}$, we shall denote by
   $\cS (U)$
the family of all complex-valued solutions $u \in C^{\infty} (U)$ of the heat
equation $u_t = \iD u$ on $U$.
The topology on $\cS (U)$ induced by embedding this space into $C^{\infty} (U)$
is actually equivalent to the topology of uniform convergence on compact
subsets of $U$.

If $\iS$ is an arbitrary subset of $\R^{n+1}$, we shall denote by $\cS (\iS)$
the family of germs on $\iS$ of solutions $u$ to the heat equations on some
open set (depending on $u$) containing $\iS$.
Such solutions form a complete locally convex space under the inductive limit
topology.

Functions in $\cS (\R^{n+1})$ will be called entire solutions of the heat
equation.
They have singularities at points at infinity of a suitable compactification
of $\R^{n+1}$.

Set
$$
   \iPhi (t,x)
 = \left\{ \begin{array}{lcl}
              \displaystyle
              \frac{1}{(4 \pi t)^{n/2}}
              \exp \Big( - \frac{|x|^2}{4 t} \Big),
            & \mbox{if}
            & t > 0,
\\
              0,
            & \mbox{if}
            & t \leq 0.
           \end{array}
   \right.
$$
This function is locally integrable over $\R^{n+1}$,
                 infinitely differentiable in $\R^{n+1} \setminus \{ 0 \}$,
and it satisfies
   $(\partial_t - \iD) \iPhi = \delta$
in the sense of distributions in $\R^{n+1}$,
   $\delta$ being Dirac's measure at $0 \in \R^{n+1}$.
The distribution $\iPhi$ is referred to as the standard fundamental solution
of convolution type for the heat equation on $\R^{n+1}$%
\footnote{Serge Lange, near the end of his life, asserted in a lecture at the
          Universit\'{e} de Montr\'{e}al, that the heat kernel is the most
          important object in all of mathematics.}.

Solutions of the heat equation of the form
\begin{equation}
\label{eq.rs}
   \sum_{j=0}^{J}
   \sum_{|\alpha| \leq A}
   c_{j,\alpha}\, \partial_t^j \partial_x^{\alpha} \iPhi (t-t_0,x-x_0),
\end{equation}
   where $c_{j,\alpha} \in \R$,
play the role of rational solutions with pole at the point
   $(t_0,x_0) \in \R^{n+1}$.
These are nothing but the potentials of distributions supported at $(t_0,x_0)$.

By a simple rational solution of the heat equation on $\R^{n+1}$ we mean any
finite linear combination
$$
   \sum_{\nu=1}^{N} c_\nu\, \iPhi (t-t_\nu,x-x_\nu)
$$
of the heat kernel itself with poles at points
   $(t_1,x_1), \ldots, (t_N,x_N)$.
From this point of view, the following is an analogue, for the heat operator,
of Runge's theorem on approximation by partial fractions for the Cauchy-Riemann
operator.

\begin{theorem}
\label{t.abpf}
For each compact set $K \subset \R^{n+1}$, the simple rational solutions to
the heat equation with poles outside of $K$ are dense in the space $\cS (K)$.
\end{theorem}

\begin{proof}
Let $u \in \cS (K)$.
We may suppose that $u \in C^{\infty} (\R^{n+1})$
   has compact support and
   satisfies $\partial_t u = \iD u$ in a neighbourhood $U$ of $K$.
Thus,
\begin{eqnarray*}
   u (t,x)
 & = &
   \int \iPhi (t-t',x-x') \left( \partial_t u - \iD u \right) (t',x') dt' dx'
\\
 & = &
   \int_{\supp (\partial_t u - \iD u)}
   \iPhi (t-t',x-x') \left( \partial_t u - \iD u \right) (t',x') dt' dx',
\end{eqnarray*}
where
   $dt' dx'$ denotes Lebesgue measure on $\R^{n+1}$.
For fixed $(t,x)$ away from the support of $\partial_t u - \iD u$, this is a
Riemann integral.

In order to estimate this integral by Riemann sums, let $Q$ be a compact set
containing the support of $\partial_t u - \iD u$ and disjoint from $K$,
such that $Q$ is a finite union of hypercubes whose sides are parallel to the
coordinate hyperplanes.
The integral over $\supp (\partial_t u - \iD u)$ is the same as the integral
over $Q$, and it is easier to take Riemann sums over $Q$.
If $\cP = \{ Q_1, \ldots, Q_N \}$ is a partition of $Q$ into hypercubes, we
denote by $|\cP|$ the mesh of $\cP$, that is, the maximum of the diameters of
the hypercubes of the partition $\cP$.
For each $Q_{\nu}$ of the partition $\cP$, choose a point
   $(t_\nu,x_\nu) \in Q_\nu$.
For each $(t,x) \in U \setminus Q$, denote by
   $\iS_{\cP} (t,x)$
the corresponding Riemann sum of the integral over $Q$.
Then, for any fixed $(t,x) \in U \setminus Q$, the Riemann sums converge to
$u (t,x)$, as $|\cP| \to 0$.

Each such Riemann sum has the form
$$
   \iS_{\cP} (t,x)
 = \sum_{\nu=1}^{N} c_\nu\, \iPhi (t-t_\nu,x-x_\nu),
$$
where the points $(t_\nu,x_\nu)$ are in $Q$.
This Riemann sum $\iS_{\cP} (t,x)$ is a simple rational solution to the heat
equation.
Thus, we have simple rational solutions to the heat equations $\iS_{\cP}$
which converge pointwise to $u$ on $U \setminus Q$, as $|\cP| \to 0$.
Since
$$
   (t,x;t',x')
 \mapsto \iPhi (t-t',x-x') \left( \partial_t u - \iD u \right) (t',x')
$$
is uniformly continuous of $(U \setminus Q) \times Q$, the convergence of this
Riemann sums (rational solutions) to $u$ is in fact uniform on $U \setminus Q$,
   as $|\cP| \to 0$.

It is not hard to see that not only do these rational solutions to the heat
equation converge to $u$, but also their partial derivatives converge uniformly
on $U \setminus Q$ to the corresponding partial derivatives of $u$.
Thus, we in fact have $C^{\infty}$ approximation.
\end{proof}

One sees that the proof actually goes through for solution of any differential
equation possessing a left fundamental solution smooth away from the diagonal.

The following result can be thought of as a theorem on rational approximation
with a fixed set of singularities.
It contains Theorem \ref{t.abpf} as a very particular case, however, the proof
is no longer constructive.

\begin{theorem}
\label{t.ra}
Let
   $K$ be a compact set in $\R^{n+1}$
and
   $\wp$ be a subset of $\R^{n+1} \setminus K$ with the property that any
   solution $g$ to $- \partial_t g = \iD g$ in $\R^{n+1} \setminus K$ which,
   together with its derivatives in $t$ up to order $J$ and in $x$ up to
   order $A$, vanishes on $\wp$, is zero in some layer around $K$.
Then the potentials (\ref{eq.rs}), where $(t_0,x_0) \in \wp$ and
                                         $c_{j,\alpha} \in \R$,
are dense in the space of all solutions to the heat equation on $K$.
\end{theorem}

By a layer around $K$ is meant any open set $U \setminus K$, where $U$ is a
neighbourhood of $K$ in $\R^{n+1}$.

\begin{proof}
Denote by $\cS (K)$ the space of all solutions of the heat equation on $K$.
Since the heat equation is hypoelliptic, each distribution $u$ satisfying
(\ref{eq.heat}) on an open set $U \subset \R^{n+1}$ belongs actually to
$C^{\infty} (U)$.
Hence $\cS (K)$ can be specified as a closed subspace of $C^{\infty} (K)$,
this latter is the space of all $C^{\infty}$ functions on neighbourhoods of
$K$ endowed with the inductive limit topology.

Let $\cR$ stand for the subspace of $\cS (K)$ consisting of the potentials
(\ref{eq.rs}), where $(t_0,x_0) \in \wp$ and
                     $c_{j,\alpha} \in \R$.
Our task is to show that $\cR$ is dense in $\cS (K)$.

To this end we use the Hahn-Banach theorem.
Pick a continuous linear functional $\mathcal{F}$ on $C^{\infty} (K)$.
Since the dual space for $C^{\infty} (K)$ just amounts to $\cE'_K$,
   the space of all distributions on $\R^{n+1}$ with support in $K$,
there is a $v \in \cE'_K$ such that
   $\mathcal{F} (u) = \langle v, u \rangle$
for all $u \in C^{\infty} (K)$.
Consider the convolution of distributions
\begin{eqnarray*}
   g (t',x')
 & = &
   \langle \iPhi (t-t',x-x'), v (t,x) \rangle
\\
 & = &
   \iPhi' \ast v,
\end{eqnarray*}
where $\iPhi' (t,x) = \iPhi (-t,-x)$ is the transposed kernel.

Since $\iPhi'$ is the fundamental solution of the transposed heat equations on
$\R^{n+1}$, it follows that
\begin{eqnarray*}
   (- \partial_{t'} - \iD_{x'}) g
 & = &
   (- \partial_{t'} - \iD_{x'}) \iPhi' \ast v
\\
 & = &
   \delta \ast v
\\
 & = &
   v
\end{eqnarray*}
on all of $\R^{n+1}$.
In particular, $g$ is a solution of the transposed heat equation in the
complement of $K$.

If $\mathcal{F}$ vanishes on $\cR$, then
\begin{eqnarray*}
   \partial_{t'}^j \partial_{x'}^{\alpha} g (t',x')
 & = &
   \langle \partial_{t'}^j \partial_{x'}^{\alpha} \iPhi (t-t',x-x'), v (t,x)
   \rangle
\\
 & = &
   (-1)^{j+|\alpha|}
   \mathcal{F}
   \left( \partial_{t}^j \partial_x^{\alpha} \iPhi (\cdot-t',\cdot-x') \right)
\\
 & = &
   0
\end{eqnarray*}
for all $(t',x') \in \wp$ and for all $j$ and
                                      $\alpha$
satisfying $j \leq J$ and
           $|\alpha| \leq A$.
By assumption, there is a neighbourhood $U$ of $K$ in $\R^{n+1}$, such that
$g = 0$ in $U \setminus K$.

Thus, if $u$ is a solution of the heat equation in a neighbourhood of $K$,
then we obtain
\begin{eqnarray*}
   \mathcal{F} (u)
 & = &
   \langle (- \partial_{t} - \iD) g, u \rangle
\\
 & = &
   \langle g, (\partial_{t} - \iD) u \rangle
\\
 & = &
   0.
\end{eqnarray*}
Now the assertion follows from the Hahn-Banach theorem.
\end{proof}

\section{Choice of pole sets}
\label{s.choice}

To effectively use Theorem \ref{t.ra} we should be able to show explicit sets
$\wp$ for which the hypotheses of the theorem are fulfilled.
In appropriate sense $\wp \in \R^{n+1} \setminus K$ are uniqueness sets for
solutions of the transposed heat equation in the complement of $K$.
More precisely,
   each solution $g$ to $- \partial_t g = \iD g$ in $\R^{n+1} \setminus K$
   satisfying $\partial_t^j \partial_x^{\alpha} g = 0$ on $\wp$ for all $j$
                                                                and $\alpha$
   with $j \leq J$ and
        $|\alpha| \leq A$
must vanish in some layer around $K$.
Such results can be derived from the analyticity of solutions to the heat
equation on characteristics.

We first recall an idea of \cite{LO} which is known as elliptic continuation
of solutions of a parabolic equation.

\begin{theorem}
\label{t.LO}
Let $g$ be a bounded solution to $\partial_t g + \iD g = 0$ in a cylinder
$Z = (0,t_0] \times B (0,R)$.
Then for each $R' < R$ there are positive constants $\varepsilon$ and
                                                    $C$
depending on $t_0$,
             $R$ and
             $R'$,
such that in 
   $Z' = (t_0-\varepsilon,t_0+\varepsilon) \times B (0,R')$
there is a solution $g'$ to the elliptic equation
   $\partial_t^2 g' + \iD g' = 0$
with Cauchy data
$$
   \begin{array}{rcl}
     g' (t_0,x)
   & =
   & g (t_0,x),
\\
     \partial_t g' (t_0,x)
   & =
   & 0
   \end{array}
$$
for $x \in B (0,R')$, satisfying
$$
   \sup_{(t,x) \in Z'} |g' (t,x)|
 \leq C\, \sup_{(t,x) \in Z} |g (t,x)|.
$$
\end{theorem}

\begin{proof}
The existence follows from the Cauchy-Kovalevskaya theorem once we observe
that the restriction $g$ to the characteristic hyperplane $t = t_0$ is real
analytic in $x \in B (0,R')$.
The estimate is a consequence of an explicit construction.
\end{proof}

From Theorem \ref{t.LO} we deduce in particular that if
   $g$ is a bounded solution to $\partial_t g + \iD g = 0$ in the cylinder $Z$
and,
   for some $x_0 \in B (0,R)$, the restriction $g (t_0,x)$ decreases faster
   than any power $|x-x_0|^k$ as $x \to x_0$,
then $g (t_0,x) \equiv 0$ for all $x \in B (0,R)$.

\begin{example}
\label{e.open}
{\em
Set $\wp = U \setminus K$ where $U$ is an arbitrary neighbourhood of $K$ in
$\R^{n+1}$.
Then the hypotheses of Theorem \ref{t.ra} are obviously fulfilled with $J = 0$
                                                                   and $A = 0$.
}
\end{example}

\begin{example}
\label{e.ruled}
{\em
Let
   $U$ be an arbitrary neighbourhood of $K$ in $\R^{n+1}$
and
   $\wp$ a subset of $U \setminus K$, such that,
   for any hyperplane $H$ orthogonal to the time axis,
   each connected component of $H \cap (U \setminus K)$ contains at least one
   point of $\wp$.
By the above, the hypotheses of Theorem \ref{t.ra} are fulfilled with $J = 0$
                                                             and $A = \infty$.
}
\end{example}

\begin{example}
\label{e.simple}
{\em
Suppose $K$ is a compact set, such that each hyperplane orthogonal to the time
axis meets $K$ in a set with connected complement.
Let $\wp$ be any curve in $\R^{n+1} \setminus K$ whose projection to the time
axis contains that of a neighbourhood of $K$.
Then the hypotheses of Theorem \ref{t.ra} are fulfilled with $J = 0$ and
                                                             $A = \infty$.
}
\end{example}

For a fuller discussion of rational approximations for elliptic equations we
refer the reader to \cite[5.3.3]{T}.

\section{Runge pairs}
\label{s.Rp}

The previous section is concerned with approximation on compact sets.
We now turn to approximation on open sets.

Let $\iO_1$ and
    $\iO_2$
be open subsets of the complex plane $\C$.
From Runge's theorem on rational approximation it follows that a necessary and
sufficient condition,
   in order that each function holomorphic in $\iO_1$ can be approximated
   uniformly on compact subsets of $\iO_1$ by functions holomorphic in $\iO_2$,
is that
   the complement $\C \setminus \iO_1$ has no compact components in $\iO_2$.
This `Runge' theorem was extended by Lax and Malgrange to approximation by
solutions of elliptic equations.

For the heat equation,
   which is of course not elliptic,
the analogous approximation problem was investigated by Jones \cite{J} and
                                                        Diaz \cite{D}.

If $\iO$ is an open set in $\R^{n+1}$, we say that $\iO$ is a Runge open set
for the heat operator, if the entire solutions to the heat equation are dense
in the solutions to the heat equation on $\iO$.
One has the following characterisation of Runge open set, proved in \cite{J}.

\begin{theorem}
\label{t.Jones}
An open set $\iO \subset \R^{n+1}$ is a Runge open set for the heat operator
if and only if,
   for every hyperplane $H$ in $\R^{n+1}$ orthogonal to the time axis,
$H \setminus \iO$ has no compact components.
\end{theorem}

Let us say that two open sets $\iO_1$ and
                              $\iO_2$
in $\R^{n+1}$,
   with $\iO_1 \subset \iO_2$,
form a Runge pair for the heat operator, if $\cS (\iO_2)$ is dense in
                                            $\cS (\iO_1)$
in the $C^{\infty}$ topology.
That is, if for each $u \in \cS (\iO_1)$ there is a sequence $\{ u_j \}$ in
                           $\cS (\iO_2)$,
such that the restrictions of the $u_j$ to $\iO_1$ and
          their partial derivatives of all orders
converge to $u$ and
            the corresponding partial derivatives of $u$,
respectively, uniformly on compact subsets of $\iO_1$.

The following theorem of Diaz \cite{D} gives a necessary and
                                               sufficient
condition that a pair $(\iO_1,\iO_2)$ of open sets in $\R^{n+1}$,
   with $\iO_1 \subset \iO_2$,
be a Runge pair for the heat operator.

\begin{theorem}
\label{t.Diaz}
Let $\iO_1 \subset \iO_2$ be open sets in $\R^{n+1}$.
A necessary and sufficient condition in order that $(\iO_1,\iO_2)$ be a Runge
pair for the heat operator is that,
   for every hyperplane $H$ in $\R^{n+1}$ orthogonal to the time axis,
the complement of $\iO_1$ in $H$ has no compact components in $\iO_2 \cap H$.
\end{theorem}

There is a small mistake in the proof of sufficiency in \cite{D}.
The assertion on p.~645, line 19, that
   an arbitrary component of $D_i \setminus Z$ has the form $A \times (a,b)$,
is not correct.
Perhaps this is easily fixed, however, we do not see how to do this.

We denote by $H (t)$ the hyperplane in $\R^{n+1}$ which orthogonal to the time
axis at time $t$.
Moreover,
   for a set $S$ in $\R^{n+1}$,
we denote by $S (t)$ the slice $S \cap H (t)$.

\begin{example}
\label{e.Diaz}
{\em
Let
   $\iO_1$ be the slit plane $\R^{2} \setminus \{ (0,x) : |x| \leq 1 \}$
and
   $\iO_2$ be the punctured plane $\R^{2} \setminus \{ (0,0) \}$.
Then, for $t \neq 0$, the set $H (t) \setminus \iO_1$ is empty, and
                              $H (0) \setminus \iO_1$
has a single compact component,
   the slit $\{ (0,x) : |x| \leq 1 \}$,
which does not belong to the set
   $\iO_2 \cap H (0) = \{ (0,x) : x \neq 0 \}$.
Hence it follows that
   the pair $(\iO_1,\iO_2)$ satisfies the condition of Theorem \ref{t.Diaz}.
However, we do not know whether or not $(\iO_1,\iO_2)$ is a Runge pair for the
heat operator in $\R^2$.
}
\end{example}

In order to highlight the problem, let us discuss some steps towards the proof
of sufficiency in Theorem \ref{t.Diaz}.
They develop Theorem 3.4.3 of \cite{Ho} which concerns general scalar
differential operators with constant coefficients in $\R^{n+1}$.
Suppose that,
   for each hyperplane $H$ orthogonal to the time axis,
every compact component of $H \setminus \iO_1$ contains
    a component of $H \setminus \iO_2$.
Choose a continuous linear functional $\cF$ on $C^{\infty} (\iO_{1})$ which
vanishes on $\cS (\iO_{2})$.
There is a distribution $v \in \cE' (\R^{n+1})$ with a compact support $K$ in
$\iO_{1}$, such that
   $\cF (u) = \langle v, u \rangle$
for all $u \in C^{\infty} (\iO_{1})$.
Since $v$ is orthogonal to all exponential solutions of the heat equation, we
conclude that there exists a distribution $g \in \cE' (\R^{n+1})$ satisfying
   $(- \partial_t - \iD) g = v$
on $\R^{n+1}$.
Now, the transposed operator $- \partial_t - \iD$ is hypoelliptic and
     $g$ satisfies $(- \partial_t - \iD) g = 0$ away from $K$,
so $g$ is an infinitely differentiable function on $\R^{n+1} \setminus K$.
Moreover, the transposed kernel $\iPhi'$ satisfies
   $\iPhi' (- \partial_t - \iD) = I$
on $\cE' (\R^{n+1})$, whence $g = \iPhi' \ast v$ on all of $\R^{n+1}$.
Since
   $(\partial_t - \iD)
    (\partial_{t'}^{j} \partial_{x'}^{\alpha} \iPhi (t-t',x-x')) = 0$
for all $(t,x) \in \iO_2$ and
        $(t',x') \in \R^2 \setminus \iO_{2}$
and all multi-indices $j$ and
                      $\alpha$
and
   $\cF$ vanishes on $\cS (\iO_{2})$,
we get
\begin{eqnarray*}
   \partial_{t'}^{j} \partial_{x'}^{\alpha} g (t',x')
 & = &
   \langle v (t,x), \partial_{t'}^{j} \partial_{x'}^{\alpha} \iPhi (t-t',x-x')
   \rangle
\\
 & = &
   0
\end{eqnarray*}
for all $(t',x')$ away from $\iO_{2}$ and
    all multi-indices $j$ and
                      $\alpha$.
Thus, $g$ vanishes to infinite order on $\R^{n+1} \setminus \iO_2$.
If we prove that $\supp g \subset \iO_{1}$, then
                 $\cF$ vanishes on $\cS (\iO_1)$
and the assertion readily follows by the Hahn-Banach theorem.
For any $t$, we denote by $\hat{K} (t)$ the topological hull of
   $K (t)$ in
   $\iO_1 (t)$.
That is,
   $\hat{K} (t)$
is the union of $K (t)$ and all components of
                $\iO_1 (t) \setminus K (t)$
which are relatively compact in $\iO_1 (t)$.
It is always a compact subset of $\iO_1 (t)$.
Since $H (t) \setminus \iO_1 (t)$ has no compact components in $\iO_2 (t)$,
the topological hull of $K (t)$ in
                        $\iO_1 (t)$
just amounts to the topological hull of $K (t)$ in
                                        $\iO_2 (t)$.
For each $t$ the function $g (t,x)$ is an analytic function of $x$ for
   $x \in H (t) \setminus K (t)$
(for solutions of the transposed heat equation are analytic in the space
 variables)
which vanishes to infinite order on $H (t) \setminus \iO_{2} (t)$ and for
all sufficiently large $x$
(because $g$ has compact support).
Denote by $\hat{K}$ the union of $\hat{K} (t)$ over all times.
By the very construction, $\hat{K}$ is a subset of $\iO_1$.
If $(t,x) \not\in \hat{K}$, then either $x$ lies in the unbounded component of
   $H (t) \setminus \hat{K} (t)$
or $x$ lies in the bounded component of
   $H (t) \setminus \hat{K} (t)$
which meets $H (t) \setminus \iO_1 (t)$ and hence
            $H (t) \setminus \iO_2 (t)$.
In the former case $g (t,x) = 0$, for the support of $g$ is compact.
In the latter case $g (t,x) = 0$, for $g (t,x)$ vanishes to infinite order on
$H (t) \setminus \iO_2 (t)$.
We thus conclude that the function $g$ vanishes away from $\hat{K}$.
Unfortunately, $\hat{K}$ fails to be a compact subset of $\iO_1$, for it is
not bounded away from the boundary of $\iO_1$.
The arguments of \cite{D} include approximation of $g$ by functions with
compact support in $\iO_1$.

\begin{remark}
\label{r.Jones}
If the slices $\iO_1 (t)$ have no ``holes,'' then,
   for each compact set $K \subset \iO_1$,
the hull $\hat{K}$ is a compact subset of $\iO_1$ as well.
Hence, the proof of sufficiency in Theorem \ref{t.Jones} is straightforward.
\end{remark}

Note that Theorem \ref{t.Diaz} is closely related to the Runge theorem for
harmonic functions, for any harmonic function in a domain $\cX \subset \R^n$
is a solution of the heat equation in the tube domain $\R \times \cX$ in
$\R^{n+1}$, which is independent of $t$.

\section{Universality}
\label{s.universal}

We shall investigate two types of universality with respect to the heat
operator, universal series and universal functions.
For an
survey of universality, we refer the reader to \cite{G}.

Let
   $V$ be a topological vector space (real or complex) and
   $\{ u_j \}_{j \in \N}$ a sequence of vectors of $V$.
The series
$$
   \sum_{j \in \N} u_j
$$
is said to be universal in $V$ if the set formed by the partial sums is dense
in $V$.

A subset $\iS$ of a topological space of Baire category II is called residual
if its complement is of Baire category I.
In this topological sense, a residual subset of a second Baire category space
contains `most' points of this space.

\begin{theorem}
\label{t.NP}
Assume that $\{ u_j \}_{j \in \N}$ is a sequence in a metrisable topological
vector space $V$ over $\mathbb{F} = \R$ or
                      $\mathbb{F} = \C$.
The following statements are equivalent:

1)
There exists a sequence $c = \{ c_j \}_{j \in \N}$ of scalars, such that the
series
$
   \displaystyle
   \sum_{j \in \N} c_j u_j
$
is universal in $V$.

2)
The set $\iS$ of such sequences $c$ is residual in $\mathbb{F}^\N$ endowed
with the product topology.

3)
For each $n \! \in \! \N$, the set of finite linear combinations from
   $\{ u_n, u_{n+1}, \ldots \}$
is dense in $V$.
\end{theorem}

\begin{proof}
See \cite{NP}.
\end{proof}

We remark that the space
$$
   \mathbb{F}^\N = \prod_{j \in \N} \mathbb{F},
$$
endowed with the product topology, is a complete metric space and hence, by
Baire's theorem, of second category.

Combining Theorem \ref{t.NP} with
          Theorem \ref{t.ra},
we arrive at universal series in spaces of solutions to the heat equation,
whose terms are simple rational solutions.

\begin{theorem}
\label{t.usors}
Suppose that
   $\iO$ is an open set in $\R^{n+1}$,
   $U$ an arbitrary neighbourhood of $\overline{\iO}$,
and
 $\{ (t_j,x_j) \}_{j \in \N}$ a dense sequence in $U \setminus \overline{\iO}$.
Then
   there is a universal series for $\cS (\iO)$, whose terms are simple rational
   solutions with poles at $\{ (t_j,x_j) \}_{j \in \N}$.
That is, there exists a sequence $\{ c_j \}_{j \in \N}$ of real numbers, such
that,
   for each $u \in \cS (\iO)$,
there is a sequence $\{ j_N \}_{N \in \N}$, such that
$$
   \sum_{j=1}^{j_N} c_j \iPhi (t-t_j,x-x_j) \to u
$$
in the topology of $C^{\infty} (\iO)$.
\end{theorem}

\begin{proof}
According to Theorem \ref{t.ra},
   for each $n \in \N$,
the set of finite linear combinations of functions
   $\{ \iPhi (t-t_n,x-x_n), \iPhi (t-t_{n+1},x-x_{n+1}), \ldots \}$
is dense in the space $\cS (\iO_1)$.
The desired conclusion now follows directly from Theorem~\ref{t.NP}.
\end{proof}

Using Examples \ref{e.ruled} and
               \ref{e.simple}
yields universal series for $\cS (\iO)$, whose terms are rational solutions
with special sets of poles outside of $\iO$.

\begin{theorem}
\label{t.us}
Suppose that
   $(\iO_1,\iO_2)$ is a Runge pair for the heat operator in $\R^{n+1}$.
Then there is a universal series for $\cS (\iO_1)$, whose terms are elements
                                  of $\cS (\iO_2)$.
More precisely,
   there exists a sequence $\{ u_j \}_{j \in \N}$ in $\cS (\iO_2)$ with the
   property that, for each solution $u \in \cS (\iO_1)$, there is a sequence
   $\{ j_N \}_{N \in \N}$, such that
$$
   \sum_{j=1}^{j_N} u_j \to u
$$
in the topology of $C^{\infty} (\iO_1)$.
\end{theorem}

\begin{proof}
Since $(\iO_1,\iO_2)$ is a Runge pair, the subspace $\cS (\iO_2)$ is dense in
                                                    $\cS (\iO_1)$.
Since
   $C^{\infty} (\iO_1)$ is separable and
   every subspace of a separable space is also separable,
we may choose a sequence $\{ u'_j \}_{j \in \N}$ in $\cS (\iO_2)$,
                                  which is dense in $\cS (\iO_1)$.
Then, for each $n \in \N$, the sequence
   $\{ u'_n, u'_{n+1}, \ldots \}$
is also dense in $\cS (\iO_1)$.
By Theorem \ref{t.NP},
   there is a sequence $\{ c_j \}_{j \in \N}$ of real numbers,
such that the series $\sum c_j u'_j$ is universal in $\cS (\iO_1)$.
To complete the proof, we set $u_j = c_j u'_j$.
\end{proof}

In order to prove the existence of universal solutions, we need a lemma on
so-called tangential approximation on an unbounded set.
In the case of the Cauchy-Riemann operator and the Laplace operator the theory
of uniform approximation has been developed not only on compact sets, but also
on (possibly unbounded) closed sets.
For the heat operator, we shall confine ourselves to very simple closed sets.
Given a sequence $\{ K_j \}$ of disjoint compacta in $\R^{n+1}$, we write
   $K_j \nearrow \infty$
if
\begin{eqnarray*}
   \sup_{(t,x) \in K_j} |(t,x)|
 & < &
   \inf_{(t,x) \in K_{j+1}} |(t,x)|,
\\
   \inf_{(t,x) \in K_j} |(t,x)|
 & \to &
   \infty
\end{eqnarray*}
as $j \to \infty$.

\begin{lemma}
\label{l.uf}
Let
   $\{ K_j \}$ be a sequence of convex compacta in $\R^{n+1}$, such that
   $K_j \nearrow \infty$.
Then, for
   each sequence $\{ u_j \}$ with $u_j \in \cS (K_j)$ and
   each sequence $\{ \varepsilon_j \}$ with $\varepsilon_j > 0$,
there is an entire solution to the heat equation $u \in \cS (\R^{n+1})$,
   such that
$$
   \sup_{(t,x) \in K_j} |u (t,x) - u_j (t,x)| < \varepsilon_j
$$
for all $j \in \N$.
\end{lemma}

\begin{proof}
For $j = 1, 2, \ldots$, choose $R_j$ such that
$$
   \sup_{(t,x) \in K_j} |(t,x)|
 < R_j
 < \inf_{(t,x) \in K_{j+1}} |(t,x)|.
$$
Set $B_0 = \emptyset$ and, for $j = 1, 2, \ldots$, let
    $B_j$ be the closed ball of radius $R_j$ centered at the origin in
    $\R^{n+1}$.
We may assume that the sequence $\{ \varepsilon_j \}$ is decreasing.
Let us proceed by induction.
Since $K_1$ is convex, by Theorem \ref{t.Jones} there exists
   $h_1 \in \cS (\R^{n+1})$
such that
   $|h_1 - u_1| < \varepsilon_1/2$ on $K_1$.
Suppose functions $h_1, \ldots, h_j$ in $\cS (\R^{n+1})$ have been constructed
with the following properties:
\begin{equation}
\label{eq.induction}
   \begin{array}{rclcl}
     \displaystyle
     | \sum_{i=1}^{j} h_i - u_j |
   & <
   & \displaystyle
     \frac{\varepsilon_j}{2^j}
   & \mbox{on}
   & K_j,
\\
     \displaystyle
     | h_j |
   & <
   & \displaystyle
     \frac{\varepsilon_j}{2^j}
   & \mbox{on}
   & B_{j-1}.
   \end{array}
\end{equation}
Since $K_{j+1}$ and
      $B_{j}$
are convex and disjoint, by Theorem \ref{t.Jones} there exists
   $h_{j+1} \in \cS (\R^{n+1})$
such that
$$
   \begin{array}{rclcl}
     \displaystyle
     | \sum_{i=1}^{j+1} h_i - u_{j+1} |
   & <
   & \displaystyle
     \frac{\varepsilon_{j+1}}{2^{j+1}}
   & \mbox{on}
   & K_{j+1},
\\
     \displaystyle
     | h_{j+1} |
   & <
   & \displaystyle
     \frac{\varepsilon_{j+1}}{2^{j+1}}
   & \mbox{on}
   & B_{j}.
   \end{array}
$$
Hence, for each $j = 1, 2, \ldots$, there exists a solution to the heat
equation
   $h_{j} \in \cS (\R^{n+1})$
which satisfies (\ref{eq.induction}).
It follows that the series
$$
   \sum_{j=1}^{\infty} h_j
$$
converges uniformly on compact subsets of $\R^{n+1}$ to an entire solution to
the heat equation
   $u \in \cS (\R^{n+1})$
with the desired properties.
\end{proof}

We are now in a position to construct a universal entire solution to the heat
equation.

\begin{theorem}
\label{t.uf}
There is an entire universal solution to the heat equation, that is, an
entire solution $u$, whose translates are dense in the space of all entire
solutions to the heat equation.
Thus,
   for each $f \in \cS (\R^{n+1})$,
there is a sequence $\{ a_j \}$ in $\R^{n+1}$, such that
   $u (\cdot + a_j) \to f (\cdot)$
in the topology of $C^{\infty} (\R^{n+1})$.
\end{theorem}

\begin{proof}
Let
   $\{ h_j \}_{j \in \N}$
be a sequence of entire solutions to the heat equation, which is dense in the
space $\cS (\R^{n+1})$ of all entire solutions to the heat equations.
Choose a sequence
   $\{ K_j \}_{j \in \N}$
of pairwise disjoint closed balls in $\R^{n+1}$, such that
   $K_j \nearrow \infty$ and
   whose radii $R_j$ tend to infinity.
Denote by $a_j$ the center of $K_j$ and set
   $u_j (t,x) = h_j ((t,x)-a_j)$.
By Lemma \ref{l.uf}, there is a solution $u \in \cS (\R^{n+1})$ such that
$$
   \sup_{(t,x) \in K_j} |u (t,x) - u_j (t,x)| < \frac{1}{j}
$$
for all $j = 1, 2, \ldots$.
Equivalently,
$$
   \sup_{|(t,x)| \leq R_j} |u ((t,x)+a_j) - h_j (t,x)| < \frac{1}{j}
$$
for all $j = 1, 2, \ldots$.
Since $1/j \to 0$,
      $R_j \to \infty$ and
      the sequence $\{ h_j \}_{j \in \N}$ is dense in $\cS (\R^{n+1})$,
the proof is complete.
\end{proof}

We remark that the abstract universality Theorem \ref{t.NP} shows that the
phenomenon of universality is generic.
That is, once we know the existence of a universal series, it turns out that
`most' series are universal.
One can also show that the result on the existence of a universal function,
Theorem \ref{t.uf} below, is generic.

\section{Burgers' equation revisited}
\label{s.Burgers}

We can now return to Burgers' equation (\ref{eq.qlpde}) in $\R^{n+1}$.
We shall tacitly assume that the coefficient $a (p)$ is independent of $p$ to
have an explicit transformation $u = U (p)$.

In order to be able to apply the inverse transformation $p = U^{-1} (u)$, we
should remain in the domain of the inverse.
It is described by
$$
   \frac{a}{U_1}\, u < 1 + U_0\, \frac{a}{U_1},
$$
$U_0$ and
$U_1$ being arbitrary constants of Example \ref{e.CHT}.
We are thus led to the study of solutions to the heat equation which take
their values on a half-axis, a fixed conical set.

Another approach we follow is to allow solutions with singularities caused
through $U$.

If $f \in \cE' (\R^{n+1})$ is an arbitrary distribution with compact support
on $\R^{n+1}$, then the inhomogeneous equation
$$
   p_t = \iD p + a (p) |\nabla p|^2 + f (t,x),
$$
cf. (\ref{eq.KPZ}), possesses a potential type solution
   $p = U^{-1} (\iPhi \ast f)$
on $\R^{n+1}$, which explicitly reads
\begin{equation}
\label{eq.rational}
   p (t,x)
 = - \frac{1}{a} \log \Big( 1 - a \frac{\iPhi \ast f\, (t,x) - U_0}{U_1} \Big).
\end{equation}
This `solution' no longer makes sense as a distribution on $\R^{n+1}$ but away
from the singularities of $f$.
If $f$ is the unit mass at a point $(t_0,x_0) \in \R^{n+1}$, we call $p (t,x)$
and its constant multiples simple rational solutions to (\ref{eq.qlpde}).
More generally, by rational solutions of (\ref{eq.qlpde}) are meant solutions
of the form (\ref{eq.rational}), where
   $f$ is a distribution on $\R^{n+1}$ whose support consists of a finite
   number of points.

From what has been proved for solutions to the heat equation we readily derive
through the transformation $u = U (p)$ analogous results for solutions to
Burgers' equation.

\begin{theorem}
\label{t.absrs}
For each compact set $K \subset \R^{n+1}$, the simple rational solutions to
(\ref{eq.qlpde}) with poles outside of $K$ are dense in the set of all smooth
solutions on $K$.
\end{theorem}

Since (\ref{eq.qlpde}) is a quasilinear equation, its solutions do not survive
under addition.
Instead, we consider compositions of solutions defined by the formula
$
   p_1 \circ p_2 := U^{-1} (U (p_1) + U (p_2)),
$
cf. Section \ref{s.C-Ht}.
As usual we introduce the composition of infinite number of solutions, to be
referred to as series, by
$$
   \circ_{j=1}^{\infty} p_j
 = \lim_{N \to \infty} p_1 \circ \ldots \circ p_N.
$$

\begin{theorem}
\label{t.usfBe}
For each open set $\iO$ in $\R^{n+1}$,
   there is a universal series in the set of solutions to (\ref{eq.qlpde})
   on $\iO$, whose terms are simple rational solutions to (\ref{eq.qlpde})
   with poles outside of $\overline{\iO}$.
\end{theorem}

That is,
   there is a sequence $\{ (t_j,x_j) \}$ of points outside of $\overline{\iO}$,
and
   a sequence $\{ c_j \}$ in $\R$,
such that,
   for each solution $p$ to (\ref{eq.qlpde}) on $\iO$,
there is a sequence $\{ j_N \}_{N \in \N}$ with
$$
   \circ_{j=1}^{j_N}
   U^{-1} \left( c_j \iPhi (t-t_j,x-x_j) \right)
 \to p
$$
in the topology of $C^{\infty} (\iO)$.

\begin{theorem}
\label{t.uffBe}
There is an entire universal solution to (\ref{eq.qlpde}), that is, an entire
solution $p$, whose translates are dense in the space of all entire solutions
to (\ref{eq.qlpde}).
\end{theorem}

This universal solution is obviously given by
   $p = U^{-1} (u)$,
where
   $u$ is an entire universal solution to the heat equation whose existence
   is guaranteed by Theorem \ref{t.uf}.

\vspace{.2cm}

{\it Acknowledgements\,}
This research was done while the first author was visiting the Universit\"{a}t
Potsdam and was supported by DFG (Deutschland) and
                             NSERC (Canada).

\newpage

\bibliographystyle{amsplain}

\end{document}